\documentclass{article}
\begin{document}
\font\germ=eufm10
\def\ssl{\hbox{\germ sl}}
\def\slh{\widehat{\ssl_2}}
\makeatletter
\def\aaa{@}
\centerline{}
\vskip1cm
\centerline{\Large\bf Polyhedral Realizations of Crystal Bases and }
\centerline{\Large\bf Braid-type Isomorphisms}
\vskip5pt
\centerline{Toshiki NAKASHIMA}
\vskip4pt
\centerline{Department of Mathematics,}
\centerline{Sophia University, Tokyo 102-8554, JAPAN}
\centerline{e-mail:\,\,toshiki@mm.sophia.ac.jp}
\makeatother

\renewcommand{\labelenumi}{$($\roman{enumi}$)$}
\renewcommand{\labelenumii}{$(${\rm \alph{enumii}}$)$}
\font\germ=eufm10

\def\al{\alpha}
\def\beneme{\begin{enumerate}}
\def\beq{\begin{equation}}
\def\beqn{\begin{eqnarray}}
\def\beqnn{\begin{eqnarray*}}
\def\bigsl{{\hbox{\fontD \char'54}}}
\def\cd{\cdots}
\def\del{\delta}
\def\Del{\Delta}
\def\ei{e_i}
\def\eit{\tilde{e}_i}
\def\eneme{\end{enumerate}}
\def\ep{\epsilon}
\def\eeq{\end{equation}}
\def\eeqn{\end{eqnarray}}
\def\eeqnn{\end{eqnarray*}}
\def\fit{\tilde{f}_i}
\def\fjt{\tilde{f}_j}
\def\ft{\tilde{f}}
\def\ge{\hbox{\germ g}}
\def\gl{\hbox{\germ gl}}
\def\hom{{\hbox{Hom}}}
\def\ify{\infty}
\def\io{\iota}
\def\kp{k^{(+)}}
\def\km{k^{(-)}}
\def\llra{\relbar\joinrel\relbar\joinrel\relbar\joinrel\rightarrow}
\def\lan{\langle}
\def\lar{\longrightarrow}
\def\lm{\lambda}
\def\Lm{\Lambda}
\def\mapright#1{\smash{\mathop{\longrightarrow}\limits^{#1}}}
\def\nd{\noindent}
\def\nn{\nonumber}
\def\oint{{\cal O}_{\rm int}(\ge)}
\def\ot{\otimes}
\def\op{\oplus}
\def\opi{\ovl\pi_{\lm}}
\def\ovl{\overline}
\def\plm{\Psi^{(\lm)}_{\io}}
\def\qq{\qquad}
\def\q{\quad}
\def\qed{\hfill\framebox[2mm]{}}
\def\QQ{\hbox{\bf Q}}
\def\qi{q_i}
\def\qii{q_i^{-1}}
\def\ran{\rangle}
\def\rlm{r_{\lm}}
\def\ssl{\hbox{\germ sl}}
\def\slh{\widehat{\ssl_2}}
\def\ti{t_i}
\def\tii{t_i^{-1}}
\def\til{\tilde}
\def\tt{{\hbox{\germ{t}}}}
\def\ttt{\hbox{\germ t}}
\def\uq{U_q(\ge)}
\def\uqm{U^-_q(\ge)}
\def\uqmq{{U^-_q(\ge)}_{\bf Q}}
\def\uqpm{U^{\pm}_q(\ge)}
\def\uqq{U_{\bf Q}^-(\ge)}
\def\uqz{U^-_{\bf Z}(\ge)}
\def\util{\widetilde\uq}
\def\vep{\varepsilon}
\def\vp{\varphi}
\def\vpi{\varphi^{-1}}
\def\xii{\xi^{(i)}}
\def\Xiioi{\Xi_{\io}^{(i)}}
\def\wtil{\widetilde}
\def\what{\widehat}
\def\wpi{\widehat\pi_{\lm}}
\def\ZZ{\hbox{\bf Z}}

\renewcommand{\thesection}{\arabic{section}}
\section{Introduction}
\setcounter{equation}{0}
\renewcommand{\theequation}{\thesection.\arabic{equation}}

In \cite{K0},\cite{K1}, Kashiwara introduced the 
theory of crystal base.
He has shown that the existence of crystal base for 
the subalgebra $\uqm$ of the given quantum algebra
$\uq$ and arbitrary  integrable highest weight $\uq$-modules.

It is a fundamental problem to present a concerete 
realization of crystal bases
for given integrable highest weight modules
as explicitly as possible.
Up to the present time, there are several  kinds of realizations, e.g.,
for finite types,  some analogues of Young tableaux were introduced
(\cite{KN}) and the piece-wise linear combinatorics for type $A$ 
 were introduced in \cite{BZ2};
for affine $A$ type, the realization of crystal bases using 
Young diagram were treated in \cite{JMMO} and for general affine cases
the new tools 'perfect crystals' and 'path realization'
were invented (\cite{KMN1}\cite{KMN2}),  and for general Kac-Moody cases
in \cite{L1}\cite{L2}
 Littelmann realized the crystal base for symmetrizable Kac-Moody
Lie algebras  and in \cite{NZ}\cite{N2} the 'polyhedral realization'
was introduced. A prticular feature of the polyhedral realization
is that it has a very explicit form and 
it is not necessary to distinguish the type of the 
underlying Kac-Moody algebra.
 But this realization
makes sense under the assumption called 'ample'
(see 3.2). We found some examples which does not hold the 
assumption.  In order to avoid this difficulties, at least 
in semi-simpe cases,  we intorduce the 'braid-type 
isomorphism' (see also \cite{L3}).

In the former half of this manuscript, 
we shall review the polyhedral realizations
of crystal bases (\cite{NZ},\cite{N2}).
It can be simply understood  as the following thing;
the crystal base of given module $V(\lm)$ or subalgebra $\uqm$
is realized as the set of lattice points 
in some convex polytope or polyhedral convex cone
in the infinite dimensional vector space $\QQ^{\ify}$.
To be more precise, we prepare some ingredients here.
Let $B_i$ ($i\in I$) be the crystal  associated with 
the simple root $\al_i$, where $I$ is the finite index set of 
simple roots (see 2.2) and $B(\lm)$ be the 
crystal base of the irreducible highest weight module $V(\lm)$
($\lm\in P_+$).
Then we have the embedding of crystal (see 3.2):
$$
\Psi^{(\lm)}_{\io}:B(\lm)\hookrightarrow
\cd\ot B_{i_k}\ot \cd\ot B_{i_2}\ot B_{i_1}\ot
R_{\lm}(\cong \ZZ^{\ify}),
\eqno(*)
$$
where $\io=\cd i_k\cd i_2 i_1$ is an infinite sequence of indices in $I$.
Here note that
each $B_{i_k}$ is identified with the set of 
integers $\ZZ$ as a set and the crystal $R_{\lm}$ includes 
only single element. 
Thus, we can identify $\ZZ^{\ify}$ with the RHS of ($*$).
Under the assumption `ample', the image of $B(\lm)$ in $\ZZ^{\ify}$ 
is obtained as the set of lattice points in the convex polyhedron
defined by some system of linear inequalities (see Sect.3).
This is the reason why it is called ``polyhedral''.

As we have mentioned above, 
in the construction of this realization, 
we do not need to distinguish the types of the associated 
Lie algebras, like 'finite', 'affine', 'hyperbolic'....
We only need  the Cartan data as the Kac-Moody algebra.
Indeed, we shall see that our realization 
can be applied to arbitrary rank 2 Kac-Moody algebras.

In the latter half of the manuscript,
we shall show the existence of the
braid-type isomorphisms between the 
tensor products of crystals $B_i$'s and $B_j$'s in the cases
$\lan h_i,\al_j\ran\lan h_j,\al_i\ran=0,1,2,3$.
The explicit form of this isomorphisms is not simple, indeed,
they are expressed by several piece-wise linear functions.
We will know, however, that it is natural in the theory of crystals.
As their application, we shall show that 
if $\ge$ is semi-simple, the crystal $B(\lm)$
can be realized in the finite rank $\ZZ$-lattice 
and the rank is equal to the length of the longest 
element in the corrsponding Weyl group
(see Proposition \ref{braid}).
By the braid-type isomorphism, we can manage to treat some
'non-ample' cases (see 4.3).

\section{Crystal Bases and Crystals}
\setcounter{equation}{0}
\renewcommand{\thesection}{\arabic{section}}
\renewcommand{\theequation}{\thesection.\arabic{equation}}

\subsection{Definition of crystal bases}

We review the crystal bases for
integrable highest weight modules and the 
nilpotent subalgebra $\uqm$ 
which are our main subject of study.
All the results in this subsection related to crystal bases
are due to M.Kashiwara \cite{K1}.

Let $\ge$ be
a  symmetrizable Kac-Moody algebra over {\bf Q}
with a Cartan subalgebra
$\ttt$, a weight lattice $P \subset \ttt^*$, the set of simple roots
$\{\al_i: i\in I\} \subset \ttt^*$,
and the set of coroots $\{h_i: i\in I\} \subset \ttt$,
where $I$ is a finite index set.
Let $\lan h,\lm\ran$ be the pairing between $\ttt$ and $\ttt^*$,
and $(\al, \beta)$ be an inner product on
$\ttt^*$ such that $(\al_i,\al_i)\in 2{\bf Z}_{\geq 0}$ and
$\lan h_i,\lm\ran={{2(\al_i,\lm)}\over{(\al_i,\al_i)}}$
for $\lm\in\ttt^*$.
Let $P^*=\{h\in \ttt: \lan h,P\ran\subset\ZZ\}$ and
$P_+:=\{\lm\in P:\lan h_i,\lm\ran\in\ZZ_{\geq 0}\}$.
We call an element in $P_+$ a {\it dominant integral weight}.
The quantum algebra $\uq$
is an associative
$\QQ(q)$-algebra generated by the $e_i$, $f_i \,\, (i\in I)$,
and $q^h \,\, (h\in P^*)$
satisfying the usual relations (see e.g.,\cite{K1} or \cite{NZ}).
The algebra $\uqm$ is the subalgebra of $\uq$ generated 
by the $f_i$ $(i\in I)$.

Let $V(\lm)$ be the irreducible highest weight module of $\uq$
with the highest weight $\lm\in P_+$.
It can be defined by
$$
V(\lm):=\uq\left/ \sum_i\uq e_i+\sum_i\uq f_i^{\lan h_i,\lm\ran+1}+
\sum_{h\in P^*}\uq(q^h-q^{\lan h,\lm\ran})\right.{.}
$$

It is well-known that as a $\uqm$-module,
there is the following natural isomorphism:
$$
V(\lm) \cong  \uqm/\sum_i\uqm f_i^{\lan h_i,\lm\ran+1}.
$$
Let $\pi_{\lm}$ be a natural projection $\uqm\longrightarrow V(\lm)$ and
set $u_{\lm}:=\pi_{\lm}(1)$. This is the unique highest weight vector in 
$V(\lm)$ up to constant.
We also denote the unit $1\in \uqm$ by $u_{\ify}$, namely,
$u_{\lm}=\pi_{\lm}(u_{\ify})$.

Let $\oint$ be the category of upper-bounded
integrable modules
 (see \cite{K1}). 
This category is semi-simple and 
each simple object is isomorphic to some $V(\lm)$.
For an object in $\oint$ (resp. $\uqm$)
and any $i\in I$, we have the decomposition:
$
V=\bigoplus_n f^{(n)}_i({\rm Ker}\,e_i)
$
(resp. $\uqm=\bigoplus_n f^{(n)}_i({\rm Ker}\,e'_i)$)
(as for $e'_i$ see \cite{K1},\cite{NZ}).
Using this, we can define the endomorphisms
$\eit$ and $\fit\in {\rm End}(V)$ 
(resp. End$(\uqm)$)
by
$$
\eit(f^{(n)}_iu)=f^{(n-1)}_iu,\,\, {\rm and}\,\,
\fit(f^{(n)}_iu)=f^{(n+1)}_iu\q {\rm for}\,\, u\in\,{\rm Ker }\,\,e_i
\,\,({\rm resp.\,\,Ker}\,\,e'_i),
$$
where we understand that $\eit u=0$ for $u\in{\rm Ker}\,e_i$
(resp. Ker $e'_i$).
Let $A\subset \QQ(q)$ be the subring of rational functions
regular at $q=0$.

In the following definition, let $M$ be 
an object in $\oint$ or $\uqm$.
\newtheorem{df2}{Definition}[section]
\begin{df2}[\cite{K1}]
A pair $(L,B)$ is a crystal base of $M$, 
if it satisfyies:
\begin{enumerate}
\item
$L$ is a free $A$-submodule of $M$
and $M\cong \QQ(q)\ot_A L$.
\item 
$L=\oplus_{\lm\in P}L_{\lm}$ and 
$B=\sqcup_{\lm\in P}B_{\lm}$ where
$L_{\lm}:=L\cap M_{\lm}$ and 
$B_{\lm}:=B\cap L_{\lm}/qL_{\lm}$.
\item
$B$ is a basis of the $\QQ$-vector space $L/qL$.
\item
$\eit L\subset L$ and
$\fit L\subset L$.

By (iv) the $\eit$ and the $\fit$ act on
$L/qL$ and
\item
$\eit B\subset B\sqcup\{0\}$ and
$\fit B\subset B\sqcup\{0\}$.
\item
For $u,v\in B$,
$\fit u=v$ if and only if  $\eit v=u$.
\end{enumerate}
\end{df2}

We set
\begin{eqnarray}
L(\lm) & :=& \sum_{i_j\in I,l\geq 0}
A\til f_{i_l}\cd \til f_{i_1}u_{\lm},
\label{def-llm}
\\
B(\lm) & := &
\{\til f_{i_l}\cd \til f_{i_1}u_{\lm}\,\,
{\rm mod}\,\,qL(\lm)\,|\,i_j\in I,l\geq 0\}
\setminus \{0\}.
\label{def-blm}
\end{eqnarray}
Here the definition of $L(\ify)$ and $B(\ify)$ are given by replacing
$u_{\lm}$ by $u_{\ify}$ in (\ref{def-llm}) and (\ref{def-blm}).

\newtheorem{thm2}[df2]{Theorem}
\begin{thm2}[\cite{K1}]
The pair $(L(\lm),B(\lm))$ (resp. $(L(\ify),B(\ify))$)
is the crystal base of $V(\lm)$ (resp. $\uqm$).
\end{thm2}

\subsection{Definition of crystals}

A crystal is a combinatorial object
obtained by abstracting the properties of crystal bases.
In what follows we fix a finite index set $I$ and a
weight lattice $P$ as above.

\begin{df2}
\label{def2.1}
A {\it crystal} $B$ is a set endowed with the following maps:
\begin{eqnarray}
&& wt:B\lar P,\\
&&\vep_i:B\lar\ZZ\sqcup\{-\infty\},\q
  \vp_i:B\lar\ZZ\sqcup\{-\infty\} \q{\hbox{for}}\q i\in I,\\
&&\eit:B\sqcup\{0\}\lar B\sqcup\{0\},
\q\fit:B\sqcup\{0\}\lar B\sqcup\{0\}\q{\hbox{for}}\q i\in I.
\end{eqnarray}
Here $0$ is an
ideal element which is not included in $B$.
Indeed, $B$ is originally a basis of a linear space,
which does not include the zero vector. This $0$ palys the simialr 
role to the zero vector.
These maps must satisfy the following axioms:
for all $b$,$b_1$,$b_2\in B$, we have
\begin{eqnarray}
&&\vp_i(b)=\vep_i(b)+\lan h_i,wt(b)\ran,
\label{vp=vep+wt}\\
&&wt(\eit b)=wt(b)+\al_i{\hbox{ if }}\eit b\in B,
\label{+alpha}\\
&&wt(\fit b)=wt(b)-\al_i{\hbox{ if }}\fit b\in B,\\
&&\eit b_2=b_1 {\hbox{ if and only if }} \fit b_1=b_2,
\label{eeff}\\
&&{\hbox{if }}\vep_i(b)=-\infty,
  {\hbox{ then }}\eit b=\fit b=0,\\
&&\eit(0)=\fit(0)=0.
\end{eqnarray}
\end{df2}

The above axioms allow us to make a crystal $B$ into
a colored oriented graph with the set of colors $I$.

\begin{df2}
\label{c-gra}
The crystal graph of a crystal $B$ is
a colored oriented graph given by
the rule : $b_1\mapright{i} b_2$ if and only if $b_2=\fit b_1$
$(b_1,b_2\in B)$.
\end{df2}

\begin{df2}
\label{df:mor}
\begin{enumerate}
\item
Let $B_1$ and $B_2$ be crystals.
A {\sl strict morphism } of crystals $\psi:B_1\lar B_2$
is a map $\psi:B_1\sqcup\{0\} \lar B_2\sqcup\{0\}$
satisfying: $\psi(0)=0$, 
if $b\in B_1$ and $\psi(b)\in B_2$, then
\begin{eqnarray}
&&\hspace{-30pt}wt(\psi(b)) = wt(b),\q \vep_i(\psi(b)) = \vep_i(b),\q
\vp_i(\psi(b)) = \vp_i(b)
\label{wt}
\end{eqnarray}
and the map $\psi: B_1\sqcup\{0\} \lar B_2\sqcup\{0\}$
commutes with all $\eit$ and $\fit$.
\item
An injective (resp. bijective)strict morphism is called an embedding 
(resp. isomorphism) of crystals.
We call $B_1$ is a subcrystal of $B_2$, if $B_1$ is a subset of $B_2$ and
becomes a crystal itself by restricting the data on it from $B_2$.
\end{enumerate}
\end{df2}

It is well-known that the algebra $\uq$ has a Hopf algebra structure.
Then the tensor product of $\uq$-modules also has
a $\uq$-module structure.
The crystal bases have very nice properties for 
tensor operations. Indeed, if $(L_i,B_i)$ is a crystal base of 
$\uq$-module $M_i$ ($i=1,2$), $(L_1\ot_A L_2, B_1\ot B_2)$
is a crystal base of $M_1\ot_{\QQ(q)} M_2$ (\cite{K1}).
Consequently, we can consider the tensor product
of crystals:
For crystals $B_1$ and $B_2$, we define their tensor product
$B_1\ot B_2$ as follows:
\begin{eqnarray}
&&B_1\ot B_2=\{b_1\ot b_2: b_1\in B_1 ,\, b_2\in B_2\},\\
&&wt(b_1\ot b_2)=wt(b_1)+wt(b_2),
\label{tensor-wt}\\
&&\vep_i(b_1\ot b_2)={\hbox{max}}(\vep_i(b_1),
  \vep_i(b_2)-\lan h_i,wt(b_1)\ran),
\label{tensor-vep}\\
&&\vp_i(b_1\ot b_2)={\hbox{max}}(\vp_i(b_2),
  \vp_i(b_1)+\lan h_i,wt(b_2)\ran),
\label{tensor-vp}\\
&&\eit(b_1\ot b_2)=
\left\{
\begin{array}{ll}
\eit b_1\ot b_2 & {\mbox{ if }}\vp_i(b_1)\geq \vep_i(b_2)\\
b_1\ot\eit b_2  & {\mbox{ if }}\vp_i(b_1)< \vep_i(b_2),
\end{array}
\right.
\label{tensor-e}
\\
&&\fit(b_1\ot b_2)=
\left\{
\begin{array}{ll}
\fit b_1\ot b_2 & {\mbox{ if }}\vp_i(b_1)>\vep_i(b_2)\\
b_1\ot\fit b_2  & {\mbox{ if }}\vp_i(b_1)\leq \vep_i(b_2).
\label{tensor-f}
\end{array}
\right.
\end{eqnarray}
Here $b_1\ot b_2$ is just another notation for an ordered pair
$(b_1,b_2)$, and we set $b_1 \ot 0 = 0 \ot b_2=0\ot 0=0$.
Note that the tensor product of crystals is
associative, namely, the crystals
$(B_1\ot B_2)\ot B_3$ and $B_1\ot(B_2\ot B_3)$ are isomorphic via
$(b_1\ot b_2)\ot b_3\leftrightarrow b_1\ot (b_2\ot b_3)$.

The examples of crystals below will be needed later.
\newtheorem{ex}[df2]{Example}
\begin{ex}
\label{Example:crystal}
\begin{enumerate}
\item
For $i\in I$, the crystal $B_i:=\{(x)_i\,: \, x \in\ZZ\}$ is defined by
\begin{eqnarray*}
&& wt((x)_i)=x \al_i,\qq \vep_i((x)_i)=-x,\qq \vp_i((x)_i)=x,\\
&& \vep_j((x)_i)=-\infty,\qq \vp_j((x)_i)
   =-\infty \q {\rm for }\q j\ne i,\\
&& \til e_j (x)_i=\del_{i,j}(x+1)_i,\qq 
\til f_j(x)_i=\del_{i,j}(x-1)_i,
\end{eqnarray*}
\item
Let $R_{\lm}:=\{r_{\lm}\}$ $(\lm\in P)$
be the crystal consisting of one-element given by (see also \cite{J}):
\begin{eqnarray*}
&& wt(r_{\lm})=\lm,\q
 \vep_i(r_{\lm})=-\lan h_i,\lm\ran,\q
   \vp_i(r_{\lm})=0,
\q \eit (r_{\lm})=\fit(r_{\lm})=0.
\end{eqnarray*}
\item
$B(\lm)$ and $B(\ify)$ can be seen as crystals by the 
following way.
We define the weight function $wt:B(\lm)\rightarrow P$
by $wt(b):=\lm-\sum_{j}\al_{i_j}$ for
$b=\til f_{i_l}\cd\til f_{i_1}u_{\lm}$ mod $qL(\lm)\ne0$.
We define integer-valued functions $\vep_i$ and $\vp_i$ on
$B(\lm)$ by
$$
\vep_i(b):={\rm max} \{k: \eit^k b\ne 0\}, \,\,
\vp_i(b):={\rm max} \{k: \fit^k b\ne 0\}.
$$

As for $B(\ify)$, the functions $wt$, $\vep_i$ and $\vp_i$ are given 
by: $wt(b):=-\sum_{j}\al_{i_j}$ 
for $b=\til f_{i_l}\cd\til f_{i_1}u_{\ify}$mod$ q L(\ify)$,
$$
\vep_i(b):={\rm max} \{k: \eit^k b\ne 0\}, \,\,
\vp_i(b):=\vep_i(b)+\lan h_i,wt(b)\ran.
$$
\end{enumerate}
\end{ex}

It is proved in \cite{K1} that
the natural projection $\pi_{\lm}:\uqm\rightarrow V(\lm)$
sends $L(\ify)$ to $L(\lm)$, and the induced map
$\what\pi_{\lm}:L(\ify)/qL(\ify)\longrightarrow L(\lm)/q L(\lm)$ sends
$B(\ify)$ to $B(\lm)\sqcup\{0\}$.
The map $\wpi$ has the following properties:
\begin{eqnarray}
&&\fit\circ\wpi=\wpi\circ\fit,
\label{fpi=pif}\\
&&\eit\circ\wpi=\wpi\circ\eit,\,\,{\rm if }\,\,\wpi(b)\ne0,
\label{epi=pie}\\
&&\hspace{-10pt}{\hbox{$\wpi:B(\ify)\setminus \{\wpi^{-1}(0)\}\longrightarrow
B(\lm)$ is bijective.}}
\label{wpi-bij}
\end{eqnarray}
Although the map $\wpi$ has such nice properties, 
it is not a strict morphism of crystals.
For instance, it does not preserve weights or 
does not necessarily commute with
the action of $\eit$ as in (\ref{epi=pie}).
We shall introduce a new strict morphism by modifying the map 
$\wpi$ in 3.2.

\section{Polyhedral Realizations of Crystal Bases}
\setcounter{equation}{0}
\renewcommand{\thesection}{\arabic{section}}
\renewcommand{\theequation}{\thesection.\arabic{equation}}

\subsection{Polyhedral Realization of $B(\ify)$}
In this subsection, we recall the results in \cite{NZ}.

We define a $\QQ(q)$-algebra anti-automorphism $*$ of $\uq$ by:
$e^*_i=e_i$,  $f^*_i=f_i$, $(q^h)^*=q^{-h}$.
This anti-automorphism has the properties (see \cite{K3}):
\begin{equation}
L(\ify)^*=L(\ify)\,\,{\rm  and}\,\,B(\ify)^*=B(\ify).
\label{*-sta}
\end{equation}
Then we can define
$\vep^*_i(b):=\vep_i(b^*)$ and $\vp^*_i(b):=\vp_i(b^*)$.

Consider the additive group
\begin{equation}
\ZZ^{\ify}
:=\{(\cd,x_k,\cd,x_2,x_1): x_k\in\ZZ
\,\,{\rm and}\,\,x_k=0\,\,{\rm for}\,\,k\gg 0\};
\label{uni-cone}
\end{equation}
we will denote by $\ZZ^{\ify}_{\geq 0} \subset \ZZ^{\ify}$
the subsemigroup of nonnegative sequences.
To the rest of this section, we fix an infinite sequence of indices
$\io=\cd,i_k,\cd,i_2,i_1$ from $I$ such that
\begin{equation}
{\hbox{
$\sharp\{k: i_k=i\}=\ify$ for any $i\in I$.}}
\label{seq-con}
\end{equation}

{\sl Remark.\,\,}
In \cite{NZ}\cite{N2}, we assume the condition $i_k\ne i_{k+1}$
for any $k$.
Note that the condition $i_k\ne i_{k+1}$ is necessarary
for removing some redundant components
(see Lemma \ref{00} below) not for the existence of the 
following embedding. 
So, without the condition $i_k\ne i_{k+1}$ 
all the results in \cite{NZ},\cite{N2} are also valid here.
We can associate to $\io$ a crystal structure
on $\ZZ^{\ify}$ and denote it by $\ZZ^{\ify}_{\io}$ 
(\cite[2.4]{NZ}).

\newtheorem{pro3}{Proposition}[section]
\begin{pro3}[\cite{K3}, See also \cite{NZ}]
\label{emb}
There is a unique embedding of crystals
\begin{equation}
\Psi_{\io}:B(\ify)\hookrightarrow \ZZ^{\ify}_{\geq 0}
\subset \ZZ^{\ify}_{\io},
\label{psi}
\end{equation}
such that
$\Psi_{\io} (u_{\ify}) = (\cd,0,\cd,0,0)$.
\end{pro3}
We call this the {\it Kashiwara embedding}
which is derived by iterating the following
type of embeddings (\cite{K3}):
\begin{enumerate}
\item
For any $i\in I$, there is a unique embedding of crystals
\begin{eqnarray}
\Psi_i : B(\infty)&\hookrightarrow &B(\infty)\ot B_i,
\label{Psii}
\end{eqnarray}
such that $\Psi_i (u_{\ify}) = u_{\ify}\ot (0)_i$.
\item
For any $b\in B(\ify)$, we can write uniquely
 $\Psi_i(b)=b'\ot \fit^m(0)_i$
where $m=\vep^*_i(b)$.
\end{enumerate}

Let us see the polyhedral realization for $B(\ify)$.
Consider the infinite dimensional vector space
$$
\QQ^{\ify}:=\{\vec{x}=
(\cd,x_k,\cd,x_2,x_1): x_k \in \QQ\,\,{\rm and }\,\,
x_k = 0\,\,{\rm for}\,\, k \gg 0\},
$$
and its dual space $(\QQ^{\ify})^*:={\rm Hom}(\QQ^{\ify},\QQ)$.
We will write a linear form $\vp \in (\QQ^{\ify})^*$ as
$\vp(\vec{x})=\sum_{k \geq 1} \vp_k x_k$ ($\vp_j\in \QQ$).

For the fixed infinite sequence
$\io=(i_k)$ we set $\kp:={\rm min}\{l:l>k\,\,{\rm and }\,\,i_k=i_l\}$ and
$\km:={\rm max}\{l:l<k\,\,{\rm and }\,\,i_k=i_l\}$ if it exists,
or $\km=0$  otherwise.
We set for $\vec x\in \QQ^{\ify}$, $\beta_0(\vec x)=0$ and
\begin{equation}
\beta_k(\vec x):=x_k+\sum_{k<j<\kp}\lan h_{i_k},\al_{i_j}\ran x_j+x_{\kp}
\qq(k\geq1).
\label{betak}
\end{equation}
We define a piecewise-linear operator $S_k=S_{k,\io}$ on $(\QQ^{\ify})^*$ by
\begin{equation}
S_k(\vp):=\left\{
\begin{array}{lll}
\vp-\vp_k\beta_k & \mbox{if}& \vp_k>0,\\
 \vp-\vp_k\beta_{\km} & \mbox{if}& \vp_k\leq 0.
\end{array}\right.
\label{Sk}
\end{equation}
Here we set
\begin{eqnarray*}
\Xi_{\io} &:=  &\{S_{j_l}\cd S_{j_2}S_{j_1}x_{j_0}\,|\,
l\geq0,j_0,j_1,\cd,j_l\geq1\},\\
\Sigma_{\io} & := &
\{\vec x\in \ZZ^{\ify}\subset \QQ^{\ify}\,|\,\vp(\vec x)\geq0\,\,{\rm for}\,\,
{\rm any}\,\,\vp\in \Xi_{\io}\}.
\end{eqnarray*}
Note that in the definition of $\Xi_{\io}$ the symbol
$x_{j_0}$ is considered as an
element in $(\QQ^{\ify})^*$, which is a function taking the
corresponding coordinate.
We impose on $\io$ the following positivity assumption:
$$
{\hbox{if $\km=0$ then $\vp_k\geq0$ for any 
$\vp\in \Xi_{\io}$ 
($\vp(\vec x)=\sum_k\vp_kx_k$)}}.
$$
\newtheorem{thm3}[pro3]{Theorem}
\begin{thm3}[\cite{NZ}]
Let $\io$ be a sequence of indices satisfying $(\ref{seq-con})$ 
and the positivity
assumption, and  $\Psi_{\io}:B(\ify)\hookrightarrow \ZZ^{\ify}_{\io}$ 
be the Kashiwara embedding associated with $\io$. Then we have 
${\rm Im}(\Psi_{\io})(\cong B(\ify))=\Sigma_{\io}$.
\end{thm3}

{\sl Remark.\,\,}
We shall see the example 
of the sequence $\io$ which does not satisfy
the positivity assumption in the end of this section.

\subsection{Polyhedral Realization of $B(\lm)$}

We review the result in \cite{N2}. 
In the rest of this section,
$\lm$ is supposed to be a dominant integral weight.
The whole story below is similar to the one of $B(\ify)$. 
The essencially different point is how to deal with the data of 
highest weight. For the purpose, the crystal $R_{\lm}$ plays 
a important role, which is 
defined in Example \ref{Example:crystal} (ii).
We shall introduce a new strict morphism of crystals by modifying 
the map $\wpi$.

Consider the crystal $B(\ify)\ot R_{\lm}$ and  define the map
\begin{equation}
\Phi_{\lm}:(B(\ify)\ot R_{\lm})\sqcup\{0\}\longrightarrow B(\lm)\sqcup\{0\},
\label{philm}
\end{equation}
by $\Phi_{\lm}(0)=0$ and $\Phi_{\lm}(b\ot r_{\lm})=\wpi(b)$ for $b\in B(\ify)$.
We set
$$
\wtil B(\lm):=
\{b\ot r_{\lm}\in B(\ify)\ot R_{\lm}\,|\,\Phi_{\lm}(b\ot r_{\lm})\ne 0\}.
$$

\begin{thm3}[\cite{N2}]
\label{ify-lm}
\begin{enumerate}
\item
The map $\Phi_{\lm}$ becomes a surjective strict morphism of crystals
$B(\ify)\ot R_{\lm}\longrightarrow B(\lm)$.
\item
$\wtil B(\lm)$ is a subcrystal of $B(\ify)\ot R_{\lm}$, 
and $\Phi_{\lm}$ induces the
isomorphism of crystals $\wtil B(\lm)\mapright{\sim} B(\lm)$.
\end{enumerate}
\end{thm3}

Let us denote $\ZZ^{\ify}_{\io}\ot R_{\lm}$ by 
$\ZZ^{\ify}_{\io}[\lm]$. Here note that
since the crystal $R_{\lm}$ has  only one element,
as a set we can identify $\ZZ^{\ify}_{\io}[\lm]$ with
$\ZZ^{\ify}_{\io}$ but their crystal structures are different.
By Theorem \ref{ify-lm}, we have the strict embedding of crystals
(see also \cite{J}):
$$
\Omega_{\lm}:B(\lm)(\cong \wtil B(\lm))\hookrightarrow B(\ify)\ot R_{\lm}.
$$
Combining $\Omega_{\lm}$ and the
Kashiwara embedding $\Psi_{\io}$,
we obtain the following:
\begin{thm3}[\cite{N2}]
\label{embedding}
There exists the unique  strict embedding of crystals
\begin{equation}
\Psi_{\io}^{(\lm)}:B(\lm)\stackrel{\Omega_{\lm}}{\hookrightarrow}
B(\ify)\ot R_{\lm}
\stackrel{\Psi_{\io}\ot {\rm id}}{\hookrightarrow}
\ZZ^{\ify}_{\io}\ot R_{\lm}=:\ZZ^{\ify}_{\io}[\lm],
\label{Psi-lm}
\end{equation}
such that $\Psi^{(\lm)}_{\io}(u_{\lm})=(\cd,0,0,0)\ot r_{\lm}$.
\end{thm3}

\vskip5pt


In the following example, we will see how the crystal 
$R_{\lm}$ works by tensoring with $B(\ify)$.
\newtheorem{ex3}[pro3]{Example}
\begin{ex3}
Let us see the simplest example $\ge=\ssl_2$-case.
Let $u_{\ify}$ be the highest weight vector in $B(\ify)$. 
Then we have 
$B(\ify)=\{\ft^n u_{\ify}\}$.
The crystal graph of $B(\ify)$ is as follows:

\nd
\unitlength 1pt
\begin{picture}(800,15)
\put(0,0){\circle{10}}
\put(0,0){\makebox(0,0){0}}
\put(5,0){\vector(1,0){30}}
\put(40,0){\circle{10}}
\put(40,0){\makebox(0,0){1}}
\put(45,0){\vector(1,0){30}}
\put(80,0){\circle{10}}
\put(80,0){\makebox(0,0){2}}
\put(85,0){\vector(1,0){30}}
\put(120,0){\circle{10}}
\put(125,0){\vector(1,0){30}}
\put(160,0){\circle{10}}
\put(165,0){\vector(1,0){30}}
\put(200,0){\circle{10}}
\put(200,0){\makebox(0,0){m}}
\put(205,0){\vector(1,0){30}}
\put(240,0){\circle{10}}
\put(245,0){\vector(1,0){30}}
\put(280,0){\circle{10}}
\put(285,0){\vector(1,0){30}}
\put(322,0){\line(1,0){3}}
\put(328,0){\line(1,0){3}}
\put(334,0){\line(1,0){3}}
\put(340,0){\line(1,0){3}}
\end{picture}
\vskip10pt

\nd
where 
\unitlength 1pt
\begin{picture}(10,5)
\put(4,4){\circle{10}}
\put(4,4){\makebox(0,0){x}}
\end{picture}$=\ft^x u_{\ify}$.

Next, let us see the crystal grah of $B(\ify)\ot R_m$ 
$(m\in\ZZ_{\geq0})$, where $R_m$ is the crystal 
as in Example \ref{Example:crystal} (ii) with
$\lm=m$.
We know that $\vp(\ft^n u_{\ify})=-n$ and $\vep(r_m)=-m$. 
Then, by (\ref{tensor-f}) we have 
$$
\ft(\ft^n u_{\ify}\ot r_m)=
\left\{
\begin{array}{lll}
\ft^{n+1}u_{\ify}\ot r_m& \mbox{if} & n<m, \cr
       \ft^n u_{\ify}\ot\ft(r_m)=0 & \mbox{if} & n \geq m.
\end{array}
\right.
$$
Thus, the crystal graph of $B(\ify)\ot R_m$ is:

\nd
\unitlength 1pt
\begin{picture}(800,15)
\put(-5,-5){\framebox(10,10){0}}
\put(5,0){\vector(1,0){30}}
\put(35,-5){\framebox(10,10){1}}
\put(45,0){\vector(1,0){30}}
\put(75,-5){\framebox(10,10){2}}
\put(85,0){\vector(1,0){30}}
\put(115,-5){\framebox(10,10){}}
\put(125,0){\vector(1,0){30}}
\put(155,-5){\framebox(10,10){}}
\put(165,0){\vector(1,0){30}}
\put(195,-5){\framebox(10,10){m}}
\put(235,-5){\framebox(10,10){}}
\put(275,-5){\framebox(10,10){}}
\put(322,0){\line(1,0){3}}
\put(328,0){\line(1,0){3}}
\put(334,0){\line(1,0){3}}
\put(340,0){\line(1,0){3}}
\end{picture}
\vskip10pt

\nd
where \framebox{x}$=\ft^x u_{\ify}\ot r_m$.
The connected component including 
$\framebox{0}=u_{\ify}\ot r_m$ is isomorphic
to the crystal $B(m)$ associated with the $m+1$-dimensional 
irreducible module $V(m)$.
The polyhedral realization below gives us the method 
how to remove the 
vectors excluded from $B(\lm)$ $($in this case, the vectors 
$\{\framebox{x}|x> m\}$.$)$.
\end{ex3}

We shall give an explicit crystal structure of
$\ZZ^{\ify}_{\io}[\lm]$ following to  \cite{N2},
which is derived by using
 the explicit formula (\ref{tensor-wt})--(\ref{tensor-f}) repeatedly.
Fix a sequence of indices $\io:=(i_k)_{k\geq 1}$ satisfying the condition
(\ref{seq-con}) and a weight $\lm\in P$.
(Here we do not necessarily assume that 
$\lm$ is dominant.)
As we stated before, 
we can identify $\ZZ^{\ify}$ with $\ZZ^{\ify}_{\io}[\lm]$
as a set. Thus $\ZZ^{\ify}_{\io}[\lm]$ can be regarded 
as a subset of $\QQ^{\ify}$, and then 
we denote an element in $\ZZ^{\ify}_{\io}[\lm]$
by $\vec x=(\cd,x_k,\cd,x_2,x_1)$.
For $\vec x=(\cd,x_k,\cd,x_2,x_1)\in \QQ^{\ify}$
we define the linear functions
\begin{eqnarray}
\sigma_k(\vec x)&:= &x_k+\sum_{j>k}\lan h_{i_k},\al_{i_j}\ran x_j
\q(k\geq1),
\label{sigma}\\
\sigma_0^{(i)}(\vec x)
&:= &-\lan h_i,\lm\ran+\sum_{j\geq1}\lan h_i,\al_{i_j}\ran x_j
\q(i\in I).
\label{sigma0}
\end{eqnarray}
Here note that
since $x_j=0$ for $j\gg0$ on $\QQ^{\ify}$,
the functions $\sigma_k$ and $\sigma^{(i)}_0$ are
well-defined.
For $\vec x\in \QQ^{\ify}$ let $\sigma^{(i)} (\vec x)
 := {\rm max}_{k: i_k = i}\sigma_k (\vec x)$, and
\begin{equation}
M^{(i)} = M^{(i)} (\vec x) :=
\{k: i_k = i, \sigma_k (\vec x) = \sigma^{(i)}(\vec x)\}.
\label{m(i)}
\end{equation}
Note that
$\sigma^{(i)} (\vec x)\geq 0$, and that
$M^{(i)} = M^{(i)} (\vec x)$ is a finite set
if and only if $\sigma^{(i)} (\vec x) > 0$.
Now we define the maps
$\eit: \ZZ^{\ify}_{\io}[\lm] \sqcup\{0\}\lar \ZZ^{\ify}_{\io}[\lm] \sqcup\{0\}$
and
$\fit: \ZZ^{\ify}_{\io}[\lm] \sqcup\{0\}\lar \ZZ^{\ify}_{\io}[\lm] \sqcup\{0\}$ 
by setting $\eit(0)=\fit(0)=0$, and 
\begin{equation}
(\fit(\vec x))_k  = x_k + \delta_{k,{\rm min}\,M^{(i)}}
\,\,{\rm if }\,\,\sigma^{(i)}(\vec x)>\sigma^{(i)}_0(\vec x);
\,\,{\rm otherwise}\,\,\fit(\vec x)=0,
\label{action-f}
\end{equation}
\begin{equation}
(\eit(\vec x))_k  = x_k - \delta_{k,{\rm max}\,M^{(i)}} \,\, {\rm if}\,\,
\sigma^{(i)} (\vec x) > 0\,\,
{\rm and}\,\,\sigma^{(i)}(\vec x)\geq\sigma^{(i)}_0(\vec x) ; \,\,
 {\rm otherwise} \,\, \eit(\vec x)=0,
\label{action-e}
\end{equation}
where $\del_{i,j}$ is the Kronecker's delta.
Here note that $\eit$ and $\fit$  act on the (infinite) tensor 
product of crystals and they act on just one component in it.
The functions $\sigma_k$ and $\sigma_0^{(i)}$ given above 
determine which component is changed by $\eit$ and $\fit$.

We also define the weight function and the functions
$\vep_i$ and $\vp_i$ on $\ZZ^{\ify}[\lm]$ by
\begin{equation}
\begin{array}{l}
wt(\vec x) :=\lm -\sum_{j=1}^{\ify} x_j \al_{i_j}, \,\,
\vep_i (\vec x) := {\rm max}(\sigma^{(i)} (\vec x),\sigma^{(i)}_0(\vec x))\\
\vp_i (\vec x) := \lan h_i, wt(\vec x) \ran + \vep_i(\vec x).
\end{array}
\label{wt-vep-vp}
\end{equation}
Now we obtain the explicit crystal structure of $\ZZ^{\ify}_{\io}[\lm]$,
Note that, in general, the subset $\ZZ^{\ify}_{\geq 0}[\lm]$ is not
a subcrystal of $\ZZ^{\ify}_{\io}[\lm]$ since it is not
stable under the action of $\eit$'s.

We fix a sequence of indices $\io$
satisfying (\ref{seq-con}) and take a dominant integral weight 
$\lm\in P_+$.
For $k\geq1$ let $k^{(\pm)}$ be  the ones in 2.4.
Let $\beta_k^{(\pm)}$ be linear functions given by
\begin{eqnarray}
&& 
\beta_k^{(+)} (\vec x)  =  \sigma_k (\vec x) - \sigma_{\kp} (\vec x)
= x_k+\sum_{k<j<\kp}\lan h_{i_k},\al_{i_j}\ran x_j+x_{\kp},
\label{beta}\\
&& \hspace{-40pt} \beta_k^{(-)} (\vec x) \nn\\
&& \hspace{-40pt}=\left\{
\begin{array}{lll}
\sigma_{\km} (\vec x) - \sigma_k (\vec x)
=x_{\km}+\sum_{\km<j<k}\lan h_{i_k},\al_{i_j}\ran x_j+x_k &
 \hspace{-10pt} \mbox{if} & \km>0,\\
\sigma_0^{(i_k)} (\vec x) - \sigma_k (\vec x)
=-\lan h_{i_k},\lm\ran+\sum_{1\leq j<k}\lan h_{i_k},\al_{i_j}\ran x_j+x_k
&
 \hspace{-10pt} \mbox{if}&  \km=0,
\end{array}
\right.
\label{beta--}
\end{eqnarray}
where the functions $\sigma_k$ and $\sigma^{(i)}_0$
are defined by (\ref{sigma}) and (\ref{sigma0}).
Here note that
$\beta_k^{(+)}=\beta_k$ and 
$\beta_k^{(-)}=\beta_{\km}  {\hbox{ \,\,if\,\, $\km>0$}}$.

Using this notation, for every $k \geq 1$, we define 
an operator
$\what S_k = \what S_{k,\io}$ for a linear function 
$\vp(\vec x)=c+\sum_{k\geq 1}\vp_kx_k$
$(c,\vp_k\in\QQ)$ on $\QQ^{\ify}$ by:

\begin{equation}
\what S_k\,(\vp) :=\left\{
\begin{array}{lll}
\vp - \vp_k \beta_k^{(+)} & \mbox{if} & \vp_k > 0,\\
\vp - \vp_k \beta_k^{(-)} & \mbox{if} & \vp_k \leq 0.
\end{array}
\right.
\label{S_k}
\end{equation}
An easy check shows that
$(\what S_k)^2=\what S_k.$

For the fixed sequence $\io=(i_k)$, 
we define $\io^{(i)}$ $(i\in I)$ to be the index $k\geq 1$
such that $i_k=i$ and $\km=0$.
Note that such $k$ is uniquely determined.
Simply,  $\io^{(i)}$ is the first
$k$ such that $i_k=i$.

Here for $\lm\in P_+$ and $i\in I$ we set
\begin{equation}
\lm^{(i)}(\vec x):=
-\beta^{(-)}_{\io^{(i)}}(\vec x)=\lan h_i,\lm\ran-\sum_{1\leq j<\io^{(i)}}
\lan h_i,\al_{i_j}\ran x_j-x_{\io^{(i)}}.
\label{lmi}
\end{equation}

For $\io$ and a dominant integral weight $\lm$,
let $\Xi_{\io}[\lm]$ be the set of all linear functions
generatd by $\what S_k=\what S_{k,\io}$ 
on the functions $x_j$ ($j\geq1$)
and $\lm^{(i)}$ ($i\in I$), namely,
\begin{equation}
\begin{array}{ll}
\Xi_{\io}[\lm]&:=\{\what S_{j_l}\cd\what S_{j_1}x_{j_0}\,
:\,l\geq0,\,j_0,\cd,j_l\geq1\}
\\
&\cup\{\what S_{j_k}\cd \what S_{j_1}\lm^{(i)}(\vec x)\,
:\,k\geq0,\,i\in I,\,j_1,\cd,j_k\geq1\}.
\end{array}
\label{Xi}
\end{equation}
Now we set
\begin{equation}
\Sigma_{\io}[\lm]
:=\{\vec x\in \ZZ^{\ify}_{\io}[\lm](\subset \QQ^{\ify})\,:\,
\vp(\vec x)\geq 0\,\,{\rm for \,\,any }\,\,\vp\in \Xi_{\io}[\lm]\}.
\label{Sigma}
\end{equation}

Here note that in general, it is possible that the set 
$\Sigma_{\io}[\lm]$ is empty. 
This never occurs in the case of $B(\ify)$.
In the case of $B(\ify)$, any element in $\Xi_{\io}$ 
has trivial constant term, but
in the base of $B(\lm)$, some element in $\Xi_{\io}[\lm]$ 
may has a negative constant term.
But here the case $\Sigma_{\io}[\lm]\ni\vec 0=(\cd,0,0)$ 
is only considered in the sequel and in this case, 
a pair
$(\io,\lm)$ is called {\it ample}.
The assumption 'ample' requires that 
the element corresponding to the highest weight vector
is contained in $\Sigma_{\io}[\lm]$.

\begin{thm3}[\cite{N2}]
\label{main}
Suppose that $(\io,\lm)$ is ample.
Let $\Psi^{(\lm)}_{\io}:B(\lm)\hookrightarrow \ZZ^{\ify}_{\io}[\lm]$
be the embedding as in (\ref{Psi-lm}). Then the image
${\rm Im}(\plm)(\cong B(\lm))$ is equal to $\Sigma_{\io}[\lm]$.
\end{thm3}

\subsection{Rank 2 case}

We apply  Theorem \ref{main} 
to the case of the
Kac-Moody algebras of rank 2.
The setting here is the same as \cite{N2}.
Without loss of generality, we can and will assume that $I=\{1,2\}$,
and $\io = (\cd,2,1,2,1)$.
The Cartan data is given by:
$$
\lan h_1,\al_1\ran= \lan h_2,\al_2\ran=2, \,\, \lan h_1,\al_2\ran=-c_1,
\,\, \lan h_2,\al_1\ran=-c_2.
$$
Here we either have $c_1 = c_2 = 0$, or both $c_1$ and $c_2$ are
positive integers.
We set $X = c_1 c_2 - 2$, and define the integer sequence
$a_l = a_l (c_1, c_2)$ for $l \geq 0$ by setting $a_0 = 0, \, a_1 = 1$
and, for $k \geq 1$,
\begin{equation}
a_{2k}  = c_1 P_{k-1} (X), \,\,
a_{2k+1} = P_k (X) + P_{k-1} (X),
\label{defcoeff}
\end{equation}
where the $P_k (X)$ are {\it Chebyshev polynomials} given 
by the following generating function:
\begin{equation}
 \sum_{k \geq 0} P_k (X) z^k = (1 - X z + z^2)^{-1}.
\label{gen-cheb}
\end{equation}
Here define $a'_l(c_1,c_2):=a_l(c_2,c_1)$.
The several first Chebyshev polynomials and terms $a_l$ are given by
$$P_0 (X) = 1, \, P_1 (X) = X, \,
P_2 (X) = X^2 - 1, \, 
P_3 (X) = X^3 - 2 X,$$
$$a_2 = c_1, \, a_3 = c_1 c_2 - 1, \, a_4 = c_1 (c_1 c_2 - 2),$$
$$a_5 = (c_1 c_2 - 1)(c_1 c_2 - 2) - 1, \, 
a_6 = c_1 (c_1 c_2 - 1)(c_1 c_2 - 3),$$
$$a_7 = c_1 c_2 (c_1 c_2 - 2)(c_1 c_2 - 3) - 1.
$$
Let $l_{\rm max} = l_{\rm max} (c_1, c_2)$ be the minimal index
$l$ such that $a_{l+1} < 0$
(if $a_l \geq 0$ for all $l \geq 0$,
then we set $l_{\rm max} = + \infty$). 
By inspection, if $c_1 c_2 = 0$ (resp. $1,2,3$) 
then $l_{\rm max} = 2$ (resp. $3, 4, 6$).
Furthermore, if $c_1 c_2 \leq 3$ then $a_{l_{\rm max}} = 0$ and 
$a_l > 0$ for $1 \leq l < l_{\rm max}$.
On the other hand, if $c_1 c_2 \geq 4$, i.e., $X \geq 2$,
it is easy to see from (\ref{gen-cheb})
that $P_k (X) > 0$ for $k \geq 0$, hence 
$a_l > 0$ for $l \geq 1$; in particular, in this case
$l_{\rm max} = + \infty$.  

\begin{thm3}[\cite{N2}]
\label{rank 2-thm}
In the rank 2 case, for a dominant integral weight 
$\lm=\lm_1\Lm_1+\lm_2\Lm_2$ $(\lm_1,\lm_2\in \ZZ_{\geq0})$ 
the image of the embedding $\Psi^{(\lm)}_{\io}$ 
is given by
\begin{equation}
{\rm Im} \,(\Psi^{(\lm)}_{\io}) = \left\{(\cd,x_2,x_1)\in\ZZ_{\geq0}^{\ify}: 
\begin{array}{l}
x_k = 0 \,\,{\rm for}\,\,k > l_{\rm max},\,\,\lm_1\geq x_1, \\
a_l x_l -a_{l-1} x_{l+1} \geq0 ,\\
\lm_2+a'_{l+1}x_l-a'_lx_{l+1}\geq0,\\
{\rm for}\,\,1 \leq l < l_{\rm max}
\end{array}
\right\}.
\label{rank 2-poly}
\end{equation}
\end{thm3}

Note that the cases when $l_{\rm max} < +\infty$, or equivalently, 
the image ${\rm Im} \,(\Psi_{\io})$ is contained in a lattice of finite rank,
just correspond to the Lie algebras 
$\ge=$ $A_1 \times A_1$, $A_2$, $B_2$ or $C_2$, $G_2$.

In conclusion of this section, we illustrate Theorem \ref{rank 2-thm}
by the example when $c_1 = c_2 = 2$, i.e., $\ge$ is the affine 
Lie algebra of type $A^{(1)}_1$.
In this case, $X = c_1 c_2 - 2 = 2$.
It follows at once from (\ref{gen-cheb}) that $P_k (2) = k+1$;
hence, (\ref{defcoeff}) gives $a_l = l$ for $l \geq 0$.
We see that for type $A^{(1)}_1$, 
$$
B(\lm)\cong 
{\rm Im} \,(\Psi^{(\lm)}_{\io}) = \{(\cd,x_2,x_1)\in\ZZ_{\geq0}^{\ify}: 
\begin{array}{l}
l x_l - (l-1) x_{l+1} \geq0, \,\,\lm_1\geq x_1\,\,{\rm and}\,\,\\
\lm_2+(l+1)x_l-lx_{l+1}\geq0 \,\,{\rm for}\,\, l \geq 1
\end{array}
\}, 
$$

\subsection{$A_n$-case}

Next, we shall apply Theorem \ref{main}
to the case when $\ge$ is of type $A_n$. 
Let us identify the index set $I$ with $[1,n] := \{1,2,\cd,n\}$ 
in the standard way; thus, the Cartan matrix 
$(a_{i,j}= \lan h_i,\al_j\ran )_{1 \leq i,j \leq n}$ is given by 
$a_{i,i}=2$, $a_{i,j}=-1$ for $|i-j|=1$, and 
$a_{i,j}=0$ otherwise. 
As the infitite sequence $\io$ let us take 
the following periodic sequence 
$$
\io = \cd,\underbrace{n,\cd,2,1}_{},
\cd,\underbrace{n,\cd,2,1}_{},\underbrace{n,\cd,2,1}_{}.
$$

Following to \cite{NZ,N2}, we shall 
change the indexing set for $\ZZ^{\ify}_{\io}$
from $\ZZ_{\geq 1}$ to $\ZZ_{\geq 1} \times [1,n]$, which is given by 
the bijection 
$\ZZ_{\geq 1} \times [1,n] \to \ZZ_{\geq 1}$ 
($(j;i) \mapsto (j-1)n + i$). 
According to this, we will write an element $\vec x \in \ZZ^{\ify}$
as a doubly-indexed family $(x_{j;i})_{j \geq 1, i \in [1,n]}$.
We will adopt the convention that $x_{j;i} = 0$ unless
$j \geq 1$ and  $i \in [1,n]$; in particular, $x_{j;0} = x_{j;n+1} = 0$
for all $j$. 

\begin{thm3} [\cite{N2}]
\label{A_n}
Let $\lm=\sum_{1\leq i\leq n}\lm_i\Lm_i$ $(\lm_i\in \ZZ_{\geq0})$ 
be a dominant integral weight. 
In the above notation, the image ${\rm Im} \,(\Psi^{(\lm)}_{\io})$ 
is the set of all integer families $(x_{j;i})$ such that 
\begin{eqnarray}
&&\hspace{-30pt}\hbox{
$x_{1;i} \geq x_{2;i-1} \geq \cd \geq x_{i;1} \geq 0$ 
for $1 \leq i \leq n$}
\label{sl-1}\\
&&\hspace{-30pt}\hbox{
$x_{j;i} = 0$ for $i+j > n+1$, }
\label{j;i=0}\\
&&\hspace{-30pt}\hbox{
$\lm_i\geq x_{j;i-j+1}-x_{j;i-j}$ for
$1\leq j\leq i\leq n$.}
\label{sl-2}
\end{eqnarray}
\end{thm3}

We give 
the example which does not satisfy the positivity assumption.

\begin{ex3}[\cite{N2}]
\label{counter-example}
We consider the case $\ge=A_3$ and take 
the sequence $\io=\cd 2\,1\,2\,3\,2\,1$,
where we do not need the explicit form of ``$\cd$'' in $\io$.
For simplicity, we write $\vec x=(\cd,x_2,x_1)$ for
an element $\vec x\in \ZZ^{\ify}_{\io}$.
In this setting, we have $\beta_1=x_1-x_2-x_4+x_5$, $\beta_2=x_2-x_3+x_4$ and 
$5^{(-)}=1$. Then 
$ S_1(x_1)=x_1-\beta_1=x_2+x_4-x_5,$ 
$S_2S_1(x_1)=x_2+x_4-x_5-\beta_2=x_3-x_5$ and 
$S_5S_2S_1(x_1)=x_3-x_5+\beta_1=x_1-x_2+x_3-x_4$.
Thus we see  $S_5S_2S_1(x_1)$ has the negative coefficient for $x_2$,
which breaks the positivity assumption.
Furthermore, this case is not ample.
Fix $\lm\in P_+$ with $\lan h_2,\lm\ran >0$. 
Since $\beta^{(-)}_2=-\lan h_2,\lm\ran+x_2-x_1$ and 
$\what S_5\what S_2\what S_1(x_1)=S_5S_2S_1(x_1)$, 
we have
$$
\what S_2\what S_5\what S_2\what S_1(x_1)=x_1-x_2+x_3-x_4+\beta^{(-)}_2
=-\lan h_2,\lm\ran+x_3-x_4,
$$
which implies 
$\vec 0=(\cd,0,0)\not\in\Sigma_{\io}[\lm]$ 
since $\lan h_2,\lm\ran>0$.
\end{ex3}

\section{Braid-type isomorphisms and its application}
\setcounter{equation}{0}
\renewcommand{\theequation}{\thesection.\arabic{equation}}
\subsection{Braid-type isomorphisms}
In this subsection we shall give the ``braid-type 
isomorphisms'' of crystals.

Let $I$ be the finite index set and  
$B_i$ and $B_j$ ($i,j\in I$) 
be the crystals as in Example \ref{Example:crystal} (i)
with the condition 
\begin{equation}
c_{ij}:=\lan h_i,\al_j\ran\lan h_j,\al_i\ran\leq 3.
\label{<3}
\end{equation}
Set $c_1:=-\lan h_i,\al_j\ran$ and $c_2:=-\lan h_j,\al_i\ran$.
In the sequel, for $x\in \QQ$ we set
$$
x_+:=\left\{
\begin{array}{lll}
x& \mbox{if} & x\geq 0,\\
0& \mbox{if} & x<0.
\end{array}
\right.
$$
\newtheorem{pr4}{Proposition}[section]
\begin{pr4}
\label{braid}
\begin{enumerate}
\item
Under the condition (\ref{<3}) there exist the following type of 
isomorphisms of crystals $\phi^{(k)}_{ij}$ $($$k=0,1,2,3$$)$:
\begin{enumerate}
\item
If $c_{ij}=0$,
$$
\phi^{(0)}_{ij}:B_i\ot B_j\mapright{\sim} B_j\ot B_i,
$$
whetre $\phi^{(0)}_{ij}((x)_i\ot(y)_j)=(y)_j\ot (x)_i$.
\item
If $c_{ij}=1$,
$$
\phi^{(1)}_{ij}:B_i\ot B_j\ot B_i\mapright{\sim} B_j\ot B_i\ot B_j,
$$
where 
$$
\phi^{(1)}_{ij}((x)_i\ot(y)_j\ot(z)_i)=
(z+(-x+y-z)_+)_j\ot (x+z)_i\ot (y-z-(-x+y-z)_+)_j.
$$
\item
If $c_{ij}=2$,
$$
\phi^{(2)}_{ij}:B_i\ot B_j\ot B_i\ot B_j
\mapright{\sim} B_j\ot B_i\ot B_j\ot B_i,
$$
where 
$\phi^{(2)}_{ij}$ is given by the following:
for $(x)_i\ot(y)_j\ot(z)_i\ot (w)_j$ we set 
$(X)_j\ot (Y)_i\ot (Z)_j\ot (W)_i
:=\phi^{(2)}_{ij}((x)_i\ot(y)_j\ot(z)_i\ot (w)_j)$.
\begin{eqnarray}
X & = & w+(-c_2x+y-w+c_2(x-c_1y+z)_+)_+,
\label{X}\\
Y & = & x+c_1w+(-x+z-c_1w+(x-c_1y+z)_+)_+,
\label{Y}\\
Z & = & y-(-c_2x+y-w+c_2(x-c_1y+z)_+)_+,
\label{Z}\\
W & = & z-c_1w-(-x+z-c_1w+(x-c_1y+z)_+)_+.
\label{W}
\end{eqnarray}
\item
If $c_{ij}=3$,  
$$
\phi^{(3)}_{ij}:B_i\ot B_j\ot B_i\ot B_j\ot B_i\ot B_j
\mapright{\sim} B_j\ot B_i\ot B_j\ot B_i\ot B_j\ot B_i,
$$
where it is defined by the following: for 
$(x)_i\ot (y)_j\ot (z)_i\ot (u)_j\ot (v)_i\ot (w)_j$
we set 
$A:=-x+c_1y-z$,
$B:=-y+c_2z-u$,
$C:=-z+c_1u-v$ and 
$D:=-u+c_2v-w$.
Then $(X)_j\ot (Y)_i\ot (Z)_j\ot (U)_i\ot (V)_j\ot (W)_i
:=\phi^{(3)}_{ij}((x)_i\ot (y)_j\ot (z)_i\ot (u)_j\ot (v)_i\ot (w)_j)$
is given by
\begin{eqnarray}
X & = & w+(D+(c_2C+(2B+c_2A_+)_+)_+)_+,\\
Y & = & x+c_1w+(c_1D+(3C+(2c_1B+2A_+)_+)_+)_+,\\
Z & = & y+u+w-X-V,\\
U & = & x+z+v-Y-W,\\
V & = & u-w-(2D+(2c_2C+(3B+c_2A_+)_+)_+)_+,\\
W & = & v-c_1w-(c_1D+(2C+(c_1B+A_+)_+)_+)_+.
\end{eqnarray}
\end{enumerate}
\item
For $k=0,1,2,3$, we have
$\phi^{(k)}_{ji}\circ\phi^{(k)}_{ij}={\rm id}$.
\end{enumerate}
\end{pr4}
\noindent
We call $\phi^{(k)}_{ij}$ a {\it braid-type isomorphism}.

In \cite{L3}, Littelmann gave the similar formula to the one in 
Proposition \ref{braid} for the 
cone realized by his path model
(he omited the complete form for the $G_2$-case.).
The different points are as follows: (1)
he obtained it by considering the actions of the Weyl groups but 
we got the formula by requiring only that it should be an 
isomorphism of crystals, and (2)
our isomorphism is for the crystal $B_i \ot B_j\ot\cd$,
which contains the cone as a subset.

\vskip10pt
{\sl Proof.}
The bijectivity of $\phi^{(k)}_{ij}$ 
follows immediately from (ii).
Thus, first let us show (ii).
The case (a) is trivial.
In the case (b), set 
$(X)_j\ot (Y)_i\ot (Z)_j
:=\phi^{(1)}_{ij}((x)_i\ot(y)_j\ot(z)_i)$ and
$(X')_i\ot (Y')_j\ot (Z')_i
:=\phi^{(1)}_{ji}((X)_j\ot(Y)_i\ot(Z)_j)$.
Then we have
$X'=Z+(-X+Y-Z)_+=y-z-(-x+y+z)_++(x-y+z)_+=y-z+(x-y+z)=x$,
$Y'=X+Z=y$ and 
$Z'=Y-Z-(-X+Y-Z)_+=x-y+2z+(-x+y-z)_+-(x-y+z)_+
=x-y+2z-(x-y+z)=z$. Here we use the formula:
$x_+-(-x)_+=x$.
Next, let us see the case (c).
We shall divide the cases by
the signs of the following four values:
$P:=x-c_1y+z$, $Q:=-c_2x+y-w$, $R:=-x+z-c_1w$
and $S:=-y+c_2z-w$. Those have the relations:
$P+c_1Q=R$, $c_2P+Q=S$, $Q+S=c_2R$ and $-P+c_1S=R$.
Set $(\al)_i\ot(\beta)_j\ot(\gamma)_j\ot (\delta)_j
=\phi^{(2)}_{ji}\circ\phi^{(2)}_{ij}((x)_i\ot(y)_j\ot(z)_i\ot(w)_j)$.
Here note that $\phi^{(2)}_{ji}$ is obtained by exchanging $c_1$ and $c_2$
in $\phi^{(2)}_{ij}$.
We shall see the case $P\leq 0$ and $Q\leq 0$.
We have $\phi^{(2)}_{ij}((x)_i\ot(y)_j\ot(z)_i\ot(w)_j))
=(w)_j\ot (x+c_1w)_i\ot (y)_i\ot(z-c_1w)_i$ and then
$\al=z-c_1w+(x-z+c_1w+c_1Q_+)_+=z-c_1w+(x-z+c_1w)_+=x$
since $x-z+c_1w=-R=-P-c_1Q\geq0$.
We also have 
$\beta=c_2x-w+(-c_2P-Q+Q_+)_+=c_2x-w+(y-c_2z+w)=y$.
Since $\al+\gamma=x+z$ and $\beta+\delta=y+w$,
we also have $\gamma=z$ and $\delta=w$.
By arguing similarly for the other cases 
we can get $\al=x$, $\beta=y$, $\gamma=z$ and $\delta=w$.

The proof for the case of $\phi^{(3)}_{ij}$ is not difficult 
but quite complicated.
It would be easier to show if we adopt another expressions
for $X,\cd,W$:
\begin{eqnarray*}
&&\hspace{-40pt}X = {\rm max}(-c_2x+y,-2y+c_2z,-c_2z+2u,-u+c_2v,w),\\
&&\hspace{-40pt}Y = {\rm max}(-x+z,x-2c_1y+3z,x-3z+2c_1u,x-c_1u+3v,x+c_1w),\\
&&\hspace{-40pt}Z = y+u+w-X-V,\\
&&\hspace{-40pt}U = x+z+v-Y-W,\\
&&\hspace{-40pt}V = {\rm min}(c_2x+w,3y-c_2z+w,2c_2z-3u+w,3u-2c_2v+w,u-w),\\
&&\hspace{-40pt}W = {\rm min}(x,c_1y-z,2z-c_1u,c_1u-2v,v-c_1w).
\end{eqnarray*}
Set
$(X')_i\ot (Y')_j\ot (Z')_i\ot (U')_j\ot (V')_i\ot (W')_j
:=\phi^{(3)}_{ji}((X)_j\ot (Y)_i\ot (Z)_j\ot (U)_i\ot (V)_j\ot (W)_i)$.
To see (ii) in this case, we may show
$(X',Y',Z',U',V',W')=(x,y,z,u,v,w)$.
For example, in the following case we shall show 
$W'=w$:
suppose that $2B+c_2A_+\geq0>3B+c_2A_+,
c_2C+2B+c_2A_+\geq0$ and $D+c_2C+2B+c_2A_+<0$.
In this case, we can show  
\begin{eqnarray}
 X=w, \,\,c_2Y-Z\geq X,\,\, 
2Z-c_2U\geq X,\,\, c_2U-2V\geq X,\,\, V-c_2W\geq X.
\label{w=w1}
\end{eqnarray}
Since $W' = {\rm min}(X,c_2Y-Z,2Z-c_2U,c_2U-2V,V-c_2W)$,
we have $W'=X=w$ by (\ref{w=w1}).
What we should show here is same as the case of $\phi^{(2)}_{ij}$,
so other cases are remained to the readers.

Let us show (i).
The bijectivity is obtained by (ii). So we shall
see that 
the map $\phi^{(k)}_{ij}$ is a strict morphism of crystals.
We should check the following:
\renewcommand{\labelenumi}{$(${\rm \arabic{enumi}}$)$}
\begin{enumerate}
\item
The map $\phi^{(k)}_{ij}$ preserves the data $wt$, $\vep_i$,
and $\vp_i$ for $i\in I$.
\item
The map $\phi^{(k)}_{ij}$ commutes all $\eit$ and $\fit$.
\end{enumerate}

Let us see (1) for the case (c).
We have $wt((x)_i\ot(y)_j\ot(z)_i\ot (w)_j)=(x+z)\al_i+(y+w)\al_j$
and $wt(\phi^{(2)}_{ij}((x)_i\ot(y)_j\ot(z)_i\ot (w)_j))=
(Y+W)\al_i+(X+Z)\al_j$.
By the explicit froms of $X,Y,Z,W$ in (\ref{X})--(\ref{W}), 
$X+Z=y+w$ and $Y+W=x+z$. Hence, we have
$wt(v)=wt(\phi^{(2)}_{ij}(v))$ for $v\in B_i\ot B_j\ot B_i\ot B_j$.
Next, let us see $\vep_i$. In the case $c_1=2$ and $c_2=1$, 
we have 
\begin{eqnarray}
\vep_i((x)_i\ot(y)_j\ot(z)_i\ot (w)_j)
&=& -x+(-x+2y-z)_+,\\
\vep_i((X)_j\ot(Y)_i\ot(Z)_j\ot (W)_i)
&=& 2X-Y+(-Y+2Z-W)_+,\\
\vep_j((x)_i\ot(y)_j\ot(z)_i\ot (w)_j)
&=& x-y+(-y+z-w)_+,\\
\vep_j((X)_j\ot(Y)_i\ot(Z)_j\ot (W)_i)
&=& -X+(-X+Y-Z)_+.
\end{eqnarray}
If $x-2y+z\leq0$, $\vep_i((x)_i\ot(y)_j\ot(z)_i\ot (w)_j)=-2x+2y-z$.
Furthermore, if $-x+y-w\geq0$, we have
$2X-Y+(-Y+2Z-W)_+=-2x+2y-x-2w-(-x+z-2w)_+
+(x-z+2w)_+=-3x+2y-2w+(x-z+2w)=-2x+2y-z$. 
If $-x+y-w\leq0$, we have
$2X-Y+(-Y+2Z-W)_+=2w-x-2w-(-x+z-2w)_+
+(-x+2y-z)_+=-2x+2y-z$
since $-x+z-2w=2(-x+y-w)+(x-2y+z)\leq0$.
Thus, we have 
$\vep_i((x)_i\ot(y)_j\ot(z)_i\ot (w)_j)
=\vep_i((X)_j\ot(Y)_i\ot(Z)_j\ot (W)_i)$ in the case 
$x-2y+z\leq0$. We can obtain other cases by the 
similar argument. 
Hence, the function $\vep_i$ is preserved by $\phi^{(2)}_{ij}$.
The case of $\vp_i$ is trivial due to the 
formula $\vp_i(b)=\lan h_i,wt(b)\ran+\vep_i(b)$,
and the fact that
the RHS of this formula is 
presereved by $\phi^{(2)}_{ij}$.
The cases for $\vep_j$ and $\vp_j$
are obtained by the similar way.
The cases for $\vep_j$ and $\vp_j$
are obtained by the similar way.
We finished (1).

Let us show (2).
For $x,y,z,w$ let $X,Y,Z,W$ be as in (\ref{X})--(\ref{W}).
We shall show the case of $\fit$. The other cases are shown similarly.

By the formula (\ref{tensor-f}) we have
\begin{eqnarray}
&&\fit((x)_i\ot(y)_j\ot(z)_i\ot (w)_j)\nn\\
&& =\left\{
\begin{array}{lll}
(x-1)_i\ot(y)_j\ot(z)_i\ot (w)_j & \mbox{if}& x-c_1y+z>0, \\
(x)_i\ot(y)_j\ot(z-1)_i\ot (w)_j & \mbox{if}& x-c_1y+z\leq0,
\end{array}
\right.
\\
&&\fit((X)_j\ot(Y)_i\ot(Z)_j\ot (W)_i)\nn\\
&&=\left\{
\begin{array}{lll}
(X)_j\ot(Y-1)_i\ot(Z)_j\ot (W)_i & \mbox{if}& Y-c_1Z+W>0,\\
(X)_j\ot(Y)_i\ot(Z)_j\ot (W-1)_i & \mbox{if}& Y-c_1Z+W\leq0.
\end{array}
\right.
\end{eqnarray}

\renewcommand{\labelenumi}{$(${\rm \roman{enumi}}$)$}
\begin{enumerate}
\item
The case $x-c_1y+z>0$:

\nd
In this case, we have
$\fit((x)_i\ot(y)_j\ot(z)_i\ot (w)_j)
=(x-1)_i\ot(y)_j\ot(z)_i\ot (w)_j$ and,
$X=w+(-y+c_2z-w)_+$,
$Y=x+c_1w+c_1(-y+c_2z-w)_+,$
$Z=y-(-y+c_2z-w)_+$ and 
$W=z-c_1w-c_1(-y+c_2z-w)_+$.
Thus, $Y-c_1Z+W=x-c_1y+z+c_1(-y+c_2z-w)_+\geq x-c_1y+w>0$ and then
we have
$\fit((X)_j\ot(Y)_i\ot(Z)_j\ot (W)_i)=
(X)_j\ot(Y-1)_i\ot(Z)_j\ot (W)_i$.
Set $(X')_j\ot(Y')_i\ot(Z')_j\ot (W')_i
:=\phi^{(2)}_{ij}((x-1)_i\ot(y)_j\ot(z)_i\ot (w)_j)$.
Since $x-c_1y+z>0$, we have $x-c_1y+z-1\geq0$.
Thus, by the explicit form of $X,Y,Z$ and $W$ we have
$X'=X$, $Y'=Y-1$, $Z'=Z$ and $W'=W$, which implies that
$\fit\circ\phi^{(2)}_{ij}=\phi^{(2)}_{ij}\circ\fit$
in the case $x-c_1y+z>0$.
\item
The case $x-c_1y+z\leq0$:

\nd
In this case, 
$\fit((x)_i\ot(y)_j\ot(z)_i\ot (w)_j)=
(x)_i\ot(y)_j\ot(z-1)_i\ot (w)_j$.
We have 
$X=w+(-c_2x+y-w)_+$,
$Y=x+c_1w+(-x+z-c_1w)_+,$
$Z=y-(-c_2x+y-w)_+$ and
$W=z-c_1w-(-x+z-c_1w)_+$.
Set $(X'')_j\ot(Y'')_i\ot(Z'')_j\ot (W'')_i
:=\phi^{(2)}_{ij}((x)_i\ot(y)_j\ot(z-1)_i\ot (w)_j)$.
If $-x+z-c_1w>0$, 
$-c_2x+y-w=\{(-x+z-c_1w)+(-x+c_1y-z)\}/c_1>0$ and then
we have
\begin{equation}
X=-c_2x+y, \q Y=z, \q Z=c_2x+w, \q W=x.
\label{XYZW}
\end{equation}
Thus, $Y-c_1Z+W=-x+z-c_1w>0$ and then
$\fit((X)_j\ot(Y)_i\ot(Z)_j\ot (W)_i)
=(X)_j\ot(Y-1)_i\ot(Z)_j\ot (W)_i$.
In this case, since $x-c_1y+z-1<0$ and $-x+(z-1)-c_1w\geq0$,
by (\ref{XYZW}) we have $X''=X$, $Y''=Y-1$, $Z''=Z$ and $W''=W$,
which implies that 
$\fit\circ\phi^{(2)}_{ij}=\phi^{(2)}_{ij}\circ\fit$
in this case.

If $-x+z-c_1w\leq0$, we  have
\begin{equation} 
X=w+(-c_2x+y-w)_+,\,
Y=x+c_1w,\q
Z=y-(-c_2x+y-w),\,
W=z-c_1w.
\label{XYZW2}
\end{equation}
Thus, we have
\begin{eqnarray*}
Y-c_1Z+W &=& x-c_1y+z+c_1(-c_2x+y-w)_+\\
&=& {\rm max}\{x-c_1y+z,-x+z-c_1w\}\leq0.
\end{eqnarray*}
So we have
$\fit((X)_j\ot(Y)_i\ot(Z)_j\ot (W)_i)
=(X)_j\ot(Y)_i\ot(Z)_j\ot (W-1)_i$.
In this case, we have still
$-x+(z-1)-c_1w<0$, 
 and then by (\ref{XYZW2}) 
$X''=X$, $Y''=Y$, $Z''=Z$ and $W''=W-1$,
which implies that 
$\fit\phi^{(2)}_{ij}=\phi^{(2)}_{ij}\fit$.
\end{enumerate}
Now we have completed to show the commutativity
of $\fit$ and $\phi^{(2)}_{ij}$.
The case of $\eit$, 
$\til e_j$ and $\til f_j$ can be shown similarly.
Now we know that the map $\phi^{(2)}_{ij}$ 
is the strict morphism of crystals and then due to 
(1) and (2), $\phi^{(2)}_{ij}$ turns out to be 
the isomorphism of crystals.

Next, we shall see that $\phi^{(3)}_{ij}$ is a strict morphism of 
crystals.
Let us show (1).
It is trivial that $\phi^{(3)}_{ij}$ preserves the
weight by its explicit form.
As for $\vep_i$, 
we have $\vep_i((x)_i\ot (y)_j\ot (z)_i\ot (u)_j\ot (v)_i\ot (w)_j)
={\mbox{max}}(-x, -x+A,-x+A+C)$.
Therefore, we consider the following three cases:
(I) $A\leq 0$ and $A+C\leq 0.$
(II) $A\geq 0 $ and $C\leq 0$.
(III)  $A\geq 0 $ and $A+C\geq 0$.
In each case, we have $\vep_i((x)_i\ot\cd\ot (w)_j)
=-x$, $-x+A$, $-x+A+C$ respectively.
And by inspection, we also have that 
$\vep_i(\phi^{(3)}_{ij}((x)_i\ot (y)_j\ot (z)_i\ot (u)_j\ot (v)_i\ot(w)_j))
=-x$, $-x+A$, $-x+A+C$ respectively,
which implies that 
$\phi^{(3)}_{ij}$ preserves $\vep_i$.
The case of $\vep_j$ can be done similarly. 
Then we also have the cases $\vp_i$ and $\vp_j$
by the formula (\ref{vp=vep+wt}).
Next, let us show (2) for $\phi^{(3)}_{ij}$.
We consider the commutativity of $\phi^{(3)}_{ij}$ and $\til f_j$
here.
Suppose that
$\fjt((x)_i\ot (y)_j\ot (z)_i\ot (u)_j\ot (v)_i\ot (w)_j)
=(x)_i\ot (y-1)_j\ot (z)_i\ot (u)_j\ot (v)_i\ot (w)_j$.
This is equivalent to the condition
$B<0$ and $B+D<0$.
We consider the cases [I] $A\leq 0$, and  [II] $A>0$.
Set
$(X')_j\ot (Y')_i\ot (Z')_j\ot (U')_i\ot (V')_j\ot (W')_i
:=\phi^{(3)}_{ij}\circ\fjt((x)_i\ot (y)_j\ot (z)_i\ot (u)_j\ot (v)_i\ot (w)_j)$.
In the case [I], we can show easily that 
$(X')_j\ot (Y')_i\ot (Z')_j\ot (U')_i\ot (V')_j\ot (W')_i
=(X)_j\ot (Y)_i\ot (Z-1)_j\ot (U)_i\ot (V)_j\ot (W)_i$
where $(X,Y,Z,U,V,W)$ is as above.
By using $A,\,B\leq0$, we can also obtain 
$\fjt\circ\phi^{(3)}_{ij}((x)_i\ot (y)_j\ot (z)_i\ot (u)_j\ot (v)_i\ot(w)_j)
=(X)_j\ot (Y)_i\ot (Z-1)_j\ot (U)_i\ot (V)_j\ot (W)_i$,
which implies $\fjt\circ\phi^{(3)}_{ij}=\fjt\circ\phi^{(3)}_{ij}$ in
this case.
In the case [II], furthermore, we consider the following three cases,
(i) $0\geq 2B+c_2A\geq 3B+c_2A$, 
(ii) $2B+c_2A\geq 0\geq 3B+c_2A$, 
(iii) $2B+c_2A\geq 3B+c_2A>0$.
In each case, we can see the commutativity of 
$\phi^{(3)}_{ij}$ and $\fjt$ by inspection.
Other cases are shown by similarly.
Now we have shown (2), and then it turns out to be that 
$\phi^{(3)}_{ij}$ is a strict morphism of crystals.
\qq\qq\qq\qq\qq\qq\qed

\subsection{Applications}
In this subsection, we introduce an application 
of the braid-type isomorphisms.

Let $\ge$ be a semi-simple Lie algebra and $W$ be the corresponding 
Weyl group. Here we denote the longest element of $W$ by
$w_0$.
Let $s_i$ ($i\in I$) be the simple reflection in $W$ and 
$N$ be the length of the longest element $w_0$ in $W$. 
Here a sequence $i_N,i_{N-1},\cd,i_2,i_1$ is called 
a {\it reduced longest word} if $s_{i_N}s_{i_{N-1}}\cd s_{i_2}s_{i_1}\in W$
is one of the reduced expressions of $w_0$.
Here we obtain
\begin{pr4}
\label{long-word}
Let $\io=i_N,i_{N-1},\cd,i_2,i_1$ $(i_j\in I)$ be one of the
reduced longest words.
Then we have
\begin{equation}
\Psi_{\io}(B(\ify))\subset
u_{\ify}\ot B_{i_N}\ot\cd\ot B_{i_2}\ot B_{i_1}
\cong \ZZ^{N},
\end{equation}
and then for $\lm\in P_+$ we also have 
\begin{equation}
\Psi^{(\lm)}_{\io}(B(\lm))\subset
u_{\ify}\ot B_{i_N}\ot\cd\ot B_{i_2}\ot B_{i_1}\ot R_{\lm}
\cong \ZZ^{N}.
\end{equation}
\end{pr4}
{\sl Remark.}\,\,
The above proposition implies that 
the crystals $B(\ify)$ and $B(\lm)$ can be embedded in the $\ZZ$-lattice 
of the 'finite' rank which is equal to the 
length of the longest element.
\vskip5pt

{\sl Proof.}\,\,\,
In order to show the proposition, 
we prepare several lemmas.
\newtheorem{lm4}[pr4]{Lemma}
\begin{lm4}
\label{00}
Let $\io=\cd,i_k,\cd,i_2,i_1$ 
be an infinite sequence of elements in $I$
such that $i_l=i_{l-1}$ for some $l>1$.
For $b\in B(\ify)$ set 
$\Psi_{\io}(b)=(\cd x_k,\cd,x_2,x_1)$.
Then we have $x_l=0$ for any $b\in B(\ify)$.
\end{lm4}

\vskip5pt

{\sl Proof.}
Let us recall the crystal structure of $\ZZ^{\ify}_{\io}$
in \cite[2.4]{NZ}.
The action of $\fit$ is determined by 
the value $\sigma_k(b)
=x_k+\sum_{j>k}\lan h_{i_k},\al_{i_j}\ran x_{i_j}$
(see (\ref{sigma}), \cite[(2.28)]{NZ}).
Then we have
$$
\sigma_{l-1}(b)-\sigma_l(b)=x_{l-1}+x_l\geq0.
$$
By (2.29) in \cite{NZ}, we know that 
any $\fit$ never acts on the $l$-th component.
Thus, we have $x_l=0$ since 
$\Psi_{\io}(b)=\til f_{i_{l}}\cd \til f_{i_{1}}(\cd,0,0)$
if $b=\til f_{i_{l}}\cd \til f_{i_{1}}u_{\ify}$.
\qed

\begin{lm4}
\label{000}
Let $\io=\cd i_k\cd i_2 i_1$ be the infinite sequence
such that
the first $N=({\rm length}\,(w_0))$ subsequence 
$i_N,\cd,i_2,i_1$ coincides with one of the reduced 
longest words.
For $b\in B(\ify)$ set $(\cd x_k\cd x_2x_1):=\Psi_{\io}(b)$.
Then we have that 
the $N+1$-th component $x_{N+1}=0$ for 
any $b\in B(\ify)$.
\end{lm4}

\vskip5pt
{\sl Proof.}
Set $i=i_{N+1}$ $(i\in I)$.
By the well-known fact (see e.g.,\cite{B}), 
there exits a reduced longest word $j_N,\cd,j_1$
such that $j_N=i$.
By applying the braid-type isomorphisms
on the components $B_{i_N}\ot\cd\ot B_{i_2}\ot B_{i_1}$
properly, we can obtain
the new tensor products of crystals
$B_{j_N}\ot \cd\ot B_{j_2}\ot B_{j_1}$
satisfying $j_N=i$.
Thus, for the new sequence 
$\io'=\cd i_{N+1},j_{N},\cd,j_2,j_1$  setting
$(\cd,x'_k,\cd,x'_2,x'_1):=\Psi_{\io'}(b)$
we have $x'_{N+1}=x_{N+1}$ since the braid-type isomorphisms
never acts on the $N+1$-th components.
Then by Lemma \ref{00}, we have $x_{N+1}=x'_{N+1}=0$.
\qed

\vskip5pt
{\sl Proof of Proposition \ref{long-word}.}
Let $\io$ and $\vec x=(\cd,x_k,\cd,x_2,x_1)$ 
be as in Lemma \ref{000}.
Let $m$-th component in $\vec x$ be the first non-trivial one 
after $N$-th component, that is 
$x_m>0$ and $x_{N+1}=\cd x_{m-1}=0$.
By Lemma \ref{000}, $m>N+1$.
Let $*$ be the map as in 3.1.
The element $b^*$ can be written uniquely:
\begin{equation}
b^*=\cd \til f_{i_k}^{x_k}\cd \til f_2^{x_2}\til f_1^{x_1}u_{\ify},
\label{**}
\end{equation}
where $\til e_{i_k}\til f_{i_{k-1}}^{x_{k-1}}
\cd \til f_2^{x_2}\til f_1^{x_1}u_{\ify}=0$
(see \cite{K3}).
Here note that 
$x_{N+1}=\cd=x_{m-1}=0$.
Adopting the new sequence 
$\io^!:=\cd i_{m} i_{N} i_{N-1}\cd i_2 i_1$
we have that the $N+1$-th component 
of $\Psi_{\io^!}(b)$ is $x_m$. 
By Lemma \ref{000}, we have $x_m=0$, which contrdicts the 
definition of $m$.
Thus, we have $x_m=0$, which implies that
all components after $N$ are just zero.
For any sequence $\cd i_k\cd i_{N+2} i_{N+1}$
the element $\cd (0)_{i_k}\ot\cd (0)_{i_{N+2}}\ot (0)_{i_{N+1}}$
can be identified with $u_{\ify}$.
Now, we have completed the proof of the 
Proposition \ref{long-word}.
\qed

\subsection{Discussions}

In Exampe \ref{counter-example}, we introduced some 
counter-example for the positivity assumption
of the polyhedral realization.
But on the other hand, 
for the sequence $\io_1=121321$, we obtain 
the image of $\Psi_{\io_1}$
by the polyhedral realization.
Here applying the braid-type isomorphism
to the case for $\io_0$ we have the one for $\io_1$, that is, 
the image of $\Psi_{\io_0}$ is given by
${\rm Im}\Psi_{\io_0}=(\phi^{(1)}_{12})_{456}{\rm Im}\Psi_{\io_1},$
where the suffix 456 menas that 
$\phi^{(1)}_{12}$ acts on the 4-th, 5-th and 6-th components
of Im$\Psi_{\io_1}$.
Indeed, we have
$$
{\rm Im}\Psi_{\io_1}
=\{(x_6,\cd,x_2,x_1)\in \ZZ^6\,|\,
x_1\geq0,\, x_2\geq x_4\geq0,\,x_3\geq x_5\geq x_6\geq0\}.
$$
Then using the explicit form of $\phi^{(1)}_{12}$ as in 
Proposition \ref{braid}, we get 
$$
{\rm Im}\Psi_{\io_0}
=\{(x_6,\cd,x_2,x_1)\in \ZZ^6\,|\,
\begin{array}{l}
x_1\geq0,\,\, x_4\geq0,\,\,x_3\geq x_4+x_6,\\
x_2+x_4\geq x_5\geq x_6\geq0,\,\,
x_2\geq x_6.
\end{array}
\}.
$$
For dominant integral weight $\lm=m_1\Lm_1+m_2\Lm_2+m_3\Lm_3$
we  also have by Theorem \ref{A_n}, 
$$
{\rm Im}\Psi_{\io_1}^{(\lm)}
=\{(x_6,\cd,x_2,x_1)\in \ZZ^6\,|\,
\begin{array}{l}
m_1\geq x_1\geq0,\, x_2\geq x_4\geq0,\,x_3\geq x_5\geq x_6\geq0,\\
m_2\geq x_2-x_1,\,x_4, m_3\geq x_3-x_2,\,x_5-x_4,\,x_6.
\end{array}
\}.
$$
Moreover, by applying $(\phi^{(1)}_{12})_{456}$ to this we have
$$
{\rm Im}\Psi_{\io_0}^{(\lm)}
=\{(x_6,\cd,x_2,x_1)\in \ZZ^6\,|\,
\begin{array}{l}
m_1\geq x_1\geq0,\, m_3\geq x_4\geq0,\,\,x_3\geq x_4+x_6,\\
x_2+x_4\geq x_5\geq x_6\geq0,\,\,x_2\geq x_6, \\
m_2\geq x_2-x_1,\,x_5-x_4,\,x_6,\,\, m_3\geq x_3-x_2.
\end{array}
\}.
$$

Now we know that the images of $\Psi_{\io_0}$ and
$\Psi_{\io_0}^{(\lm)}$
is realized 
in the polyhedral convex cone or convex polytope in $\ZZ^6$.
But in general, it is not clear 
since the braid-type isomorphisms are 
not linear but piece-wise linear.
(For A-type it seems to be correct.).
Our further problem is 
to describe explicitly the image of $\Psi_{\io}$ and $\Psi^{(\lm)}_{\io}$
for an arbitrary reduced longest word $\io$.

\end{document}